\documentclass{article}
\usepackage{amssymb}

\usepackage{graphicx}
\usepackage{amsmath}
\usepackage[symbol]{footmisc}
\usepackage[french]{babel}

\newtheorem{theorem}{Theorem}[section]

\newtheorem{conjecture}[theorem]{Conjecture}
\newtheorem{corollary}[theorem]{Corollary}

\newtheorem{lemma}[theorem]{Lemma}

\newtheorem{proposition}[theorem]{Proposition}
\newtheorem{remark}[theorem]{Remark}

\newenvironment{proof}[1][Proof]{\textbf{#1.} }{\ \rule{0.5em}{0.5em}}
\input{tcilatex}
\renewcommand{\theequation}{\thesection.\arabic{equation}}

\begin{document}

\begin{center}
{\huge On a particular specialization of monomial symmetric functions}

\bigskip

{\large Vincent Brugidou}$^{a}$%
\footnotetext{\textit{E-mail address }vincent.brugidou@univ-lille.fr}
\end{center}

\bigskip

\begin{center}
$^{a}$\textit{\ Universit\'{e} de Lille, 59655 Villeneuve d'Ascq cedex,
France}
\end{center}

\bigskip

\textbf{Abstract: }Let $m_{\lambda }$ be the monomial symmetric functions, $%
\lambda $ being an integer partition of $n\in \mathbb{N}^{\ast }$. For the
specialization corresponding to the $q$-deformation of the exponential, we
prove that each $m_{\lambda }$ is associated with a polynomial $J_{\lambda
}\left( q\right) $ whose coefficients belong to $\mathbb{Z}$. $J_{\lambda }$
is an algebraic generalization of the case $\lambda =\left( n\right) $ for
which $J_{\left( n\right) }=J_{n}$ is the enumerator of tree inversions. We
establish several relations between $J_{\lambda }$ and the polynomials $%
J_{n}^{\left( r\right) }$, which were defined algebraically in our previous
work for $n\geq r\geq 1$, and are classical combinatorial enumerators with $%
J_{n}^{\left( 1\right) }=J_{n}$. Based on inductive computations of $%
J_{\lambda }$ for $n\leq 43$, we conjecture that the coefficients of each $%
J_{\lambda }$\ are strictly positive and log-concave. As a consequence of
Huh's Theorem on the $h$-vector of matroid complex, it is shown that the
coefficients of $J_{n}^{\left( r\right) }$ are strictly positive and
log-concave, which gives a second argument in favor of these conjectures.
Furthermore, we prove that the last $n-1$ coefficients of $J_{\lambda }$ are
proportional to the first coefficients of the ($n-r-1$)-th column of
Pascal's triangle, where $r$ is the length of $\lambda $. This is a third
argument in favor of the conjectures. However, The calculation of $J_{\left(
3,2,1\right) }$ reveals a difficulty, if one wishes to prove the conjectures
through the application of Huh's theorem mentioned above.\bigskip

\textit{keywords: }monomial symmetric functions, q-deformation of the
exponential, tree inversions, log-concavity, $h$-vector of matroids.

\section{\protect\bigskip Introduction}

Let $p_{n}^{\left( r\right) }$\ be the symmetric functions defined for any
pair of integers $\left( n,r\right) $ such that $n\geq r\geq 1$, by

\begin{equation}
p_{n}^{(r)}=\sum\limits_{\left| \lambda \right| =n,\;l\left( \lambda \right)
=r}m_{\lambda }\text{,}  \tag{1.1}
\end{equation}
where $m_{\lambda \text{ }}$ are the monomial symmetric functions and the
sum is taking over all integer partitions $\lambda $ of $n$ of length $%
l(\lambda )=r$. The functions $p_{n}^{(r)}$ are introduced with this
notation in $\left[ 12\text{, Example 19, p. 33}\right] $. In particular, $%
p_{n}^{\left( 1\right) }$ is the power symetric function $p_{n}=m_{n}$. In
[8], we showed the following result for the specialization whose generating
function of elementary symmetric function is given by the following Equation 
$\left( 1.2\right) $, sometimes called the $q$-deformation of the
exponential.

\bigskip

\textbf{Theorem 6.3 i) of }$\left[ 8\right] $ \textit{For }$n\geq r\geq 1$, 
\textit{and the specialization } 
\begin{equation}
E_{xq}(t)=\sum\limits_{n=0}^{\infty }q^{\binom{n}{2}}\dfrac{t^{n}}{n!}\text{,%
}  \tag{1.2}
\end{equation}
\textit{\ we have:} 
\begin{equation}
p_{n}^{(r)}=\left( 1-q\right) ^{n-r}\dfrac{q^{\binom{r}{2}}}{r!\left(
n-r\right) !}J_{n}^{\left( r\right) }(q)\text{,}  \tag{1.3}
\end{equation}
\textit{where }$J_{n}^{\left( r\right) }$\textit{\ is a monic polynomial
with positive integer coefficients, whose constant term is }$\left(
n-r\right) !$\textit{\ \ and whose degree is }$\binom{n-1}{2}-\binom{r-1}{2}$%
\textit{. Moreover, \ for all }$r\geq 1$,\textit{\ }$J_{r}^{\left( r\right)
}=1$\textit{.\bigskip }

When $r=1$, this gives 
\begin{equation}
m_{n}=p_{n}=p_{n}^{\left( 1\right) }=\dfrac{\left( 1-q\right) ^{n-1}}{\left(
n-1\right) !}J_{n}\text{,}  \tag{1.4}
\end{equation}
where $J_{n}=J_{n}^{\left( 1\right) }$ is the enumerator polynomial of
inversions in trees on $n$ vertices, introduced in [13]. More generally, the
polynomials $J_{n}^{\left( r\right) }$ enumerate\ inversions for sequences
for sequences of\ ''colored'' forests introduced by Stanley and Yan (see $%
\left[ 18\right] $ and $\left[ 19\right] $), or level statistic enumerators
introduced in Section 9 of $\left[ 8\right] $. Alternatively, the reciprocal
polynomial of $J_{n}^{\left( r\right) }$, denoted $\overline{J_{n}^{\left(
r\right) }}$ in $\left[ 8\right] $, is the sum enumerator of generalized
parking functions. We refer to the survey [20] and to $\left[ 8\right] $ for
more details on these combinatorial representations.

When we take $\lambda $ to be the partition $\left( n\right) $, the monomial 
$m_{\lambda }$ factors as in $\left( 1.4\right) $ for the specialization $%
\left( 1.2\right) $. The question that motivated our work is wether a
similar factorization exists for any partition $\lambda $. The answer is
affirmative and is given by Theorem 3.1 and Corollary 3.5. The factorization
of $m_{\lambda }$ yields a polynomial $J_{\lambda }$ whose coefficients
belong to $\mathbb{Z}$, and which algebraically generalizes $J_{n}$. This is
proven by induction using two order relations on the set of integer
partitions. It follows two recurrence relations allowing all $J_{\lambda }$
to be calculated from the initial condition $J_{1}=1$ (Corollary 3.7). In
Section 4, based on the computations of $J_{\lambda }$ for $n=\left| \lambda
\right| \leq 43$, we conjecture that, for any partition $\lambda $, the
coefficients of $J_{\lambda }$ are positive and log-concave. A second
argument supporting these conjectures is that we show (Proposition 4.4) the
log-concavity of $J_{n}^{\left( r\right) }$ using Huh's results on the $h$%
-vector of the matroid complex of representable matroids over a field of
characteristic zero. We also prove in Section 5 that the last $n-1$
coefficients of $J_{\lambda }$ are proportional to the first $n-1$\
coefficients of the ($n-r-1$)-th column of Pascal's triangle. This provides
a third argument in favor of the conjectures since it is well known that
these columns are log-concave. To conclude, in Section 6 we discuss various
approaches to fully prove the conjectures.

\section{\protect\bigskip Prerequisites}

\QTP{Body Math}
We assume the reader has some familiarity with integer partitions and
symmetric functions, as presented in Chap. 1 of $\left[ 12\right] $. Here
are our notations. $\mathcal{P}$ and $\mathcal{P}_{n}$ are respectively the
set of all the integer partitions and the set of integer partitions of $n\in 
\mathbb{N}^{\ast }=\mathbb{N}-\left\{ 0\right\} $.If $\lambda =\left(
\lambda _{1},\lambda _{2},...,\lambda _{r}\right) $ is a partition of the
integer $n$, then, $\left| \lambda \right| =\lambda _{1}+\lambda
_{2}+...+\lambda _{r}=n$,$\;l(\lambda )=r$, $\lambda !=\lambda _{1}!\lambda
_{2}!...\lambda _{r}!$ and

\QTP{Body Math}
\begin{equation}
n(\lambda )=\sum\limits_{i\geq 1}(i-1)\lambda _{i}=\sum\limits_{i\geq 1}%
\binom{\lambda _{i}^{\prime }}{2}\text{,}  \tag{2.1}
\end{equation}
where $\lambda ^{\prime }=\left( \lambda _{1}^{\prime },\lambda _{2}^{\prime
},..\right) $ is the conjugate partition of $\lambda $. If $r_{i}$ is the
number of parts of $\lambda $ equal to $i\in \mathbb{N}^{\ast }$, the
sequence of multiplicity is $m\left( \lambda \right) =\left(
r_{1},r_{2},...\right) $, and we set $\left| m\left( \lambda \right) \right|
=r_{1}+r_{2}+...=r$ and $m\left( \lambda \right) !=r_{1}!r_{2}!...$.. We
also write $\lambda =1^{r_{1}}2^{r_{2}}...$.

\QTP{Body Math}
In $\mathcal{P}_{n}$, $\preceq $, $\leq $ and $\sqsubseteq $ will designate
the reverse lexicographic order, the dominance order and the refinement
order, respectively. It is known ($\left[ 12\right] $, Chap.1) that\ in $%
\mathcal{P}_{n}$ 
\begin{equation}
\lambda \sqsubseteq \mu \Rightarrow \lambda \leq \mu \Rightarrow \lambda
\preceq \mu \text{.}  \tag{2.2}
\end{equation}

\QTP{Body Math}
We partially generalize the above definitions to any strict composition of $%
n $, i.e. any r-multiplet $u=\left( u_{1},u_{2},...,u_{r}\right) $ of
integers strictly greater than zero such that $u_{1}+u_{2}+...+u_{r}=n$, by
setting $\left| u\right| =u_{1}+u_{2}+...+u_{r}$, $l(u)=r$. Furthermore, $%
\lambda \left( u\right) $ is defined as the partition of $n$, which is
composed of the $u_{i}$ arranged in a non-increasing way. We clearly have $%
\left| \lambda \left( u\right) \right| =\left| u\right| $ and $l\left(
\lambda \left( u\right) \right) =l\left( u\right) $.

\QTP{Body Math}
If $K$ is a commutative field, $\Lambda _{K}$ is the algebra of symmetric
functions in the indeterminates $X=\left( x_{i}\right) _{i\geq 1}$with
coefficients in $K$. For the purposes of this article $K$ will be $\mathbb{Q}
$ or $\mathbb{Q(}q)$, the rational fractions in the indeterminate $q$. For $%
\lambda \in \mathcal{P}$, $\left( m_{\lambda }\right) ,\left( e_{\lambda
}\right) $ and $\left( p_{\lambda }\right) $ are the three classical bases
of $\Lambda _{K}$. Agreeing that $e_{0}=1$,$\;$we recall that in $\Lambda
_{K}\left[ \left[ t\right] \right] $

\QTP{Body Math}
\begin{equation}
E(t)=\sum\limits_{n\geq 0}e_{n}t^{n}=\prod\limits_{i\geq 1}(1+x_{i}t)\text{.}%
\;  \tag{2.3}
\end{equation}

\QTP{Body Math}
Basic knowledge is required on partially ordered sets (posets), matroids and
the Tutte polynomial, as presented, for example, in the corresponding
Wikipedia articles. If $n\in \mathbb{N}$, $\left[ n\right]
_{q}=1+q+q^{2}+...+q^{n-1}$ is the $q$-analogue of $n$. If $P\left( q\right) 
$ is a polynomial in $q$, then $\left\langle q^{m}\right\rangle P\left(
q\right) $ denotes the coefficient of $q^{m}$ in $P\left( q\right) $, and
the order of $P$, denoted $ord\left( P\right) $, is the smallest power of $q$
with a non zero coefficient. Finally, if $E$ is a finite set, $\left|
E\right| $ denotes the cardinality of $E$.

\section{$\protect\bigskip $Expression of $m_{\protect\lambda }$ for the
specialization $e_{n}=q^{\binom{n}{2}}/n!$}

Let define the augmented monomial symmetric functions by $\ $ 
\begin{equation}
\widetilde{m_{\lambda }}=m\left( \lambda \right) !\;m_{\lambda }\text{.} 
\tag{3.1}
\end{equation}

\begin{theorem}
\bigskip For any partition $\lambda $ and for the specialization $E_{xq}(t)$%
, $\widetilde{m_{\lambda }}$ is given by 
\begin{equation}
\widetilde{m_{\lambda }}=\left( 1-q\right) ^{\left| \lambda \right| -l\left(
\lambda \right) }M_{\lambda }(q)\text{,}  \tag{3.2}
\end{equation}
where $M_{\lambda }(q)$ belongs to $\mathbb{Q}\left[ q\right] $ and
satisfies:

$a)$ The degree of $M_{\lambda }$ is 
\begin{equation}
d(\lambda )=\binom{\left| \lambda \right| -1}{2}+l\left( \lambda \right) -1%
\text{,}  \tag{3.3}
\end{equation}

and 
\begin{equation}
\left\langle q^{d\left( \lambda \right) }\right\rangle M_{\lambda }=\dfrac{%
\left( l\left( \lambda \right) -1\right) !}{\left( \left| \lambda \right|
-1\right) !}\text{.}  \tag{3.4}
\end{equation}
$\qquad b)$ The order of $M_{\lambda }$ is 
\begin{equation}
ord\left( M_{\lambda }\right) =n\left( \lambda \right) \text{,}  \tag{3.5}
\end{equation}

and 
\begin{equation}
\left\langle q^{n\left( \lambda \right) }\right\rangle M_{\lambda }=\dfrac{%
m\left( \lambda \right) !}{\lambda ^{\prime }!}\text{.}  \tag{3.6}
\end{equation}
$\qquad c)$ For $\lambda =\left( n\right) $ we have 
\begin{equation}
M_{\left( n\right) }=\dfrac{J_{n}}{\left( n-1\right) !}\text{.}  \tag{3.7}
\end{equation}
\end{theorem}

\bigskip Note that the highest degree monomial of $M_{\lambda }$ only
depends on $\left| \lambda \right| $ and $l\left( \lambda \right) $. We need
two lemmas to prove this theorem.

\begin{lemma}
The function $f(x)=\left( x-1\right) \left( x-2\right) /2$ is superadditive
on the interval $\left[ 1,+\infty \right[ $, i.e. for all $x,y$ in this
interval 
\begin{equation*}
f(x+y)\geq f\left( x\right) +f(y)\text{.}
\end{equation*}
Moreover, the inequality is strict for $x$ or $y$ greater\ than 1.
\end{lemma}

\begin{proof}
\bigskip By developing, it follows 
\begin{equation*}
f(x+y)=\dfrac{x^{2}+2xy+y^{2}-3x-3y+2}{2}=f(x)+f(y)+\left( xy-1\right) \text{%
,}
\end{equation*}
and $xy\geq 1$ for $x,y\geq 1$, the inequality being strict if $x$ or $y>1$.
\end{proof}

\begin{lemma}
\bigskip Let $\lambda $ and $\mu $ be two partitions of $n$, then: \ $\mu
<\lambda \Rightarrow n\left( \mu \right) >n\left( \lambda \right) $
\end{lemma}

\begin{proof}
By transitivity of the order relation, it suffices to show the inequality
for $\ \lambda $ covering $\mu $. Thus necessarily $\lambda =R_{i,j}\mu $
where $R_{i,j}$ is the raising operator of Mac Donald $\left[ 12\text{,
Chap. 1}\right] $. If $\mu =\left( \mu _{1},\mu _{2},...,\mu _{i},...,\mu
_{j},..,\mu _{r}\right) $, we therefore\ have $\lambda =(\mu _{1},\mu
_{2},...,\mu _{i}+1,...,\mu _{j}-1,...,\mu _{r})$, from where 
\begin{equation*}
n(\lambda )=\sum\limits_{i\geq 1}\left( i-1\right) \mu _{i}+\left(
i-1\right) -\left( j-1\right) =n\left( \mu \right) +i-j<n\left( \mu \right) 
\text{.}
\end{equation*}
\end{proof}

\begin{proof}[Proof of Theorem 3.1]
\bigskip\ 

$i)$ If $\lambda =\left( n\right) $ which corresponds to $c)$, then $r=1$
and according to $\left( 1.4\right) $: 
\begin{equation*}
\widetilde{m}_{\left( n\right) }=m_{n}=p_{n}=\left( 1-q\right) ^{n-1}\dfrac{%
J_{n}\left( q\right) }{\left( n-1\right) !}\text{,}
\end{equation*}
where $J_{n}$ is a monic polynomial of degree $\binom{n-1}{2}$, with an
order equal to 0 and a constant term equal to $\left( n-1\right) !$. Taking $%
M_{\left( n\right) }=J_{n}/\left( n-1\right) !$, Equation $\left( 3.2\right) 
$ and points $a)$ and $b)$ are verified in this case and so $c)$ also.

$ii)$ Equation $\left( 3.2\right) $\ \ and $a)$\ will be proved in the
general case by induction with the total order $\trianglelefteq $ defined on 
$\mathcal{P}$\ by 
\begin{equation*}
\lambda \vartriangleleft \mu \Leftrightarrow \left\{ 
\begin{array}{c}
\text{\ \ \ }\;\left| \lambda \right| <\left| \mu \right| \text{, } \\ 
\text{or}\;\left| \lambda \right| =\left| \mu \right| \text{ and }\mu \prec
\lambda \text{.}
\end{array}
\right. 
\end{equation*}

The first partitions, thus ordered, are : $\left( 1\right) \vartriangleleft
\left( 2\right) \vartriangleleft \left( 1,1\right) \vartriangleleft \left(
3\right) \vartriangleleft \left( 2,1\right) \vartriangleleft \left(
1,1,1\right) \vartriangleleft ...$

The case $\lambda =\left( 1\right) $ which corresponds to $i)$ with $n=1$,
has already been proven. Assume that $\left( 3.2\right) $ and $a)$ are
proved for all partition $\mu $ such that $\mu \vartriangleleft \lambda $
with $\left| \lambda \right| =n$ and let prove them for $\lambda $. If $%
\lambda =\left( n\right) $ then it is proved by case $i)$. So suppose that $%
\lambda =\left( \lambda _{1},...,\lambda _{r}\right) $, with $l\left(
\lambda \right) =r\geq 2$. By Theorem 1 of $\left[ 14\right] $, we have: 
\begin{equation}
\widetilde{m}_{\lambda }=\widetilde{m}_{\lambda _{r}}\widetilde{m}_{\lambda
^{\ast }}-\sum\limits_{i=1}^{r-1}\widetilde{m}_{\lambda ^{\left( i\right) }}
\tag{3.8}
\end{equation}
with $\lambda ^{\ast }=\left( \lambda _{1},...,\lambda _{r-1}\right) $ \ and 
$\lambda ^{\left( i\right) }=\lambda (\left( \lambda _{1},...,\lambda
_{i-1},\lambda _{i}+\lambda _{r},\lambda _{i+1},...,\lambda _{r-1}\right) )$%
. In others words, $\lambda ^{\ast }$ is the partition obtained by removing
the last part of $\lambda $, and $\lambda ^{\left( i\right) }$ is obtained
by adding to the ith part of $\lambda ^{\ast }$ the last part of $\lambda $,
and by reordering the strict composition of $n$ thus obtained in a
non-increasing way.

We have $\left| \lambda ^{\ast }\right| =\left| \lambda \right| -\lambda
_{r}<n$ and $\lambda _{r}<n$ therefore $\lambda ^{\ast }\vartriangleleft
\lambda $ and $\left( \lambda _{r}\right) \vartriangleleft \lambda $.
Similarly for $i$ between $1$ and $r-1$ we clearly have $\left| \lambda
^{\left( i\right) }\right| =n$ and $\lambda \sqsubset \lambda ^{\left(
i\right) }$ which implies with $\left( 2.2\right) $ $\lambda \prec \lambda
^{\left( i\right) }$. Since $\left| \lambda ^{\left( i\right) }\right|
=\left| \lambda \right| $ this leads with $\left( 2.2\right) $\ to $\lambda
^{\left( i\right) }\vartriangleleft \lambda .$ From $l\left( \lambda ^{\ast
}\right) =l\left( \lambda ^{\left( i\right) }\right) =r-1$, it follows by
induction hypothesis that 
\begin{equation}
\left\{ 
\begin{array}{c}
\widetilde{m}_{\lambda _{r}}=\left( 1-q\right) ^{\lambda _{r}-1}M_{\left(
\lambda _{r}\right) }(q)\text{,} \\ 
\widetilde{m}_{\lambda ^{\ast }}=\left( 1-q\right) ^{n-\lambda
_{r}-r+1}M_{\lambda ^{\ast }}(q)\text{,} \\ 
\widetilde{m}_{\lambda ^{\left( i\right) }}=\left( 1-q\right)
^{n-r+1}M_{\lambda ^{\left( i\right) }}(q)\text{,}
\end{array}
\right.   \tag{3.9}
\end{equation}
the three polynomials $M_{\left( \lambda _{r}\right) }$, $M_{\lambda ^{\ast
}}$, $M_{\lambda ^{\left( i\right) }}$ verifying the equations of $a).$ By
replacing these expressions in $\left( 3.8\right) $ we obtain, any performed
calculation : 
\begin{equation*}
\widetilde{m}_{\lambda }=\left( 1-q\right) ^{n-r}M_{\lambda }(q)\text{,}
\end{equation*}
with 
\begin{equation}
M_{\lambda }=M_{\left( \lambda _{r}\right) }M_{\lambda ^{\ast }}+\left(
q-1\right) \sum\limits_{i=1}^{r-1}M_{\lambda ^{\left( i\right) }}\text{.} 
\tag{3.10}
\end{equation}
By induction hypothesis, the degrees of $M_{\left( \lambda _{r}\right) }$, $%
M_{\lambda ^{\ast }}$, $M_{\lambda ^{\left( i\right) }}$ are $\binom{\lambda
_{r}-1}{2}$, $\binom{n-\lambda _{r}-1}{2}+r-2$ and $\binom{n-1}{2}+r-2$,
respectively. So the degree of $M_{\left( \lambda _{r}\right) }M_{\lambda
^{\ast }}$ is $\binom{\lambda _{r}-1}{2}+\binom{n-\lambda _{r}-1}{2}+r-2$.
The leading coefficient of $M_{\lambda ^{\left( i\right) }}$ being equal by
induction hypothesis to $\left( r-2\right) !/\left( n-1\right) !$, \ the
degree of $\sum\limits_{i=1}^{r-1}M_{\lambda ^{\left( i\right) }}$ is also
equal to $\binom{n-1}{2}+r-2$ and its leading coefficient is clearly $\left(
r-1\right) !/\left( n-1\right) !$. The degree of $M_{\left( \lambda
_{r}\right) }M_{\lambda ^{\ast }}$ and that of $\sum\nolimits_{i=1}^{r-1}M_{%
\lambda ^{\left( i\right) }}$ are therefore in the same order as $\binom{%
\lambda _{r}-1}{2}+\binom{n-\lambda _{r}-1}{2}$ and $\binom{n-1}{2}$. Let $%
x=\lambda _{r}$, $y=n-\lambda _{r}$ and $f(x)=\left( x-1\right) \left(
x-2\right) /2$, the degree of $M_{\left( \lambda _{r}\right) }M_{\lambda
^{\ast }}$ and that of $\sum\nolimits_{i=1}^{r-1}M_{\lambda ^{\left(
i\right) }}$ are in the same order as $f\left( x\right) +f\left( y\right) $
et $f(x+y)=f(n)$. By Lemma 3.2, the degree of $\sum\nolimits_{i=1}^{r-1}M_{%
\lambda ^{\left( i\right) }}$ is then greater than or equal to that of $%
M_{\left( \lambda _{r}\right) }M_{\lambda ^{\ast }}$. It follows that the
degree of $M_{\lambda }$ is equal to that of $\left( q-1\right)
\sum\nolimits_{i=1}^{r-1}M_{\lambda ^{\left( i\right) }}$ and is thus equal
to $\binom{n-1}{2}+r-1$ which is the value given by $\left( 3.3\right) $.
The leading coefficient of $M_{\lambda }$ is that of $\left( q-1\right) $ $%
\sum\nolimits_{i=1}^{r-1}M_{\lambda ^{\left( i\right) }}$ and is thus equal
to $\left( r-1\right) !/\left( n-1\right) !$, which is the value given by $%
\left( 3.4\right) $.

$iii)$ Let now prove $b)$ of Theorem 3.1. According to Equation $\left(
2.3\right) $ p. 20 of $\left[ 12\right] $, for any partition $\lambda $ we
have

\begin{equation}
m_{\lambda }=e_{\lambda ^{\prime }}-\sum\limits_{\mu <\lambda }a_{\lambda
\mu }m_{\mu }\text{.}  \tag{3.11}
\end{equation}
For $E_{xq}\left( q\right) $, it follows from $\left( 1.2\right) $ and $%
\left( 2.1\right) $: 
\begin{equation}
e_{\lambda ^{\prime }}=\prod\limits_{i\geq 1}e_{\lambda _{i}^{\prime
}}=\prod\limits_{i\geq 1}\dfrac{q^{\binom{\lambda _{i}^{^{\prime }}}{2}}}{%
\lambda _{i}^{\prime }!}=\dfrac{q^{n(\lambda )}}{\lambda ^{\prime }!}\text{.}
\tag{3.12}
\end{equation}
We do an induction on the order relation $\preceq $ in $\mathcal{P}_{n}$.
For $\lambda =1^{n}$, we have $\widetilde{m}_{1^{n}}=n!m_{1^{n}}=n!e_{n}=q^{%
\binom{n}{2}}$ and on the other hand, $m_{1^{n}}=\left( 1-q\right)
^{n-n}M_{1^{n}}$. It follows $M_{1^{n}}=q^{\binom{n}{2}}$, which verifies
the equations of $b).$

Let $\lambda \in \mathcal{P}_{n}$ with $\lambda \neq 1^{n}$and suppose that
the equations of $b)$ have been proven for all partition $\mu $ such that $%
\left| \mu \right| =n$ and $\mu \prec \lambda $. In the sum of $\left(
3.11\right) $, $\mu <\lambda $ implies $\mu \prec \lambda $ by $\left(
2.2\right) $. Then by induction hypothesis $ord\left( M_{\mu }\right) =n(\mu
)$, and 
\begin{equation*}
\widetilde{m}_{\mu }=\left( 1-q\right) ^{\left| \mu \right| -l(\mu )}M_{\mu
}\Rightarrow ord\left( m_{\mu }\right) =ord\left( \widetilde{m}_{\mu
}\right) =ord\left( M_{\mu }\right) =n\left( \mu \right) \text{.}
\end{equation*}
By Lemma 3.3, $\mu <\lambda \Rightarrow n\left( \mu \right) >n\left( \lambda
\right) $, therefore $\left( 3.11\right) $ shows that the monomial of
minimal degree of $m_{\lambda }$ is $e_{\lambda ^{\prime }}=q^{n\left(
\lambda \right) }/\lambda ^{\prime }!$. So we checked that 
\begin{equation*}
ord\left( M_{\lambda }\right) =ord\left( m_{\lambda }\right) =n\left(
\lambda \right) \;\text{and \ }\left\langle q^{n\left( \lambda \right)
}\right\rangle M_{\lambda }=m\left( \lambda \right) !\left\langle q^{n\left(
\lambda \right) }\right\rangle e_{\lambda ^{\prime }}=\dfrac{m\left( \lambda
\right) !}{\lambda ^{\prime }!}.
\end{equation*}
\end{proof}

\begin{corollary}
\bigskip Let $\left( n,r\right) $ be a pair of positive integers, then: 
\begin{equation}
\text{For \ }n\geq r\geq 1\text{,}\;\;\;q^{\binom{r}{2}}\dfrac{J_{n}^{\left(
r\right) }}{r!\left( n-r\right) !}=\sum\limits_{\left| \lambda \right|
=n,\;l\left( \lambda \right) =r}\dfrac{M_{\lambda }}{m\left( \lambda \right)
!}\text{.}  \tag{3.13}
\end{equation}
\begin{equation}
\text{For }n-1\geq r\geq 1\text{,}\;\ \ \ \ \ \;\;r!\sum\limits_{\left|
\lambda \right| =n,\;l\left( \lambda \right) =r}\dfrac{M_{\lambda }}{m\left(
\lambda \right) !}=q^{\binom{r}{2}}\sum\limits_{\left| \mu \right| =n-r}%
\left[ r\right] _{q}^{l\left( \mu \right) }\dfrac{M_{\mu }}{m\left( \mu
\right) !}\text{.}  \tag{3.14}
\end{equation}
In particular for $r=1$ and $n\geq 2$ 
\begin{equation}
\dfrac{J_{n}}{\left( n-1\right) !}=M_{\left( n\right) }=\sum\limits_{\left|
\lambda \right| =n-1}\dfrac{M_{\lambda }}{m\left( \lambda \right) !}. 
\tag{3.15}
\end{equation}
\end{corollary}

\begin{proof}
\bigskip $\left( 3.13\right) $ is obtained by replacing in $\left(
1.1\right) $ $p_{n}^{\left( r\right) }$ and $m_{\lambda }$ by the right side
of $\left( 1.3\right) $ and $\left( 3.2\right) $, respectively. For Eq. $%
\left( 3.14\right) $, we start from the following linear recurrence proved
in ($\left[ 8\right] $, Equation $\left( 6.5\right) $): 
\begin{equation*}
\text{For }n-1\geq r\geq 1\text{,}\;\;\;J_{n}^{\left( r\right)
}(q)=\sum\limits_{j=1}^{n-r}\left[ r\right] _{q}^{j}\;q^{\binom{j}{2}}\binom{%
n-r}{j}J_{n-r}^{\left( j\right) }(q)\text{.}
\end{equation*}
and we replace $J_{n}^{\left( r\right) }$ and $J_{n-r}^{\left( j\right) }$
with their expressions in termes of $M_{\lambda }$, coming from $\left(
3.13\right) $. With $r=1$, the first equality of $\left( 3.15\right) $\
follows from $\left( 3.13\right) $ and the second from $\left( 3.14\right) $.
\end{proof}

\begin{corollary}
\bigskip For the specialization $E_{xq}\left( t\right) $, we have the
following generalization of $\left( 1.4\right) $ for any partition $\lambda $%
: $\ $%
\begin{equation}
m_{\lambda }=\left( 1-q\right) ^{\left| \lambda \right| -l(\lambda )}\dfrac{%
q^{n\left( \lambda \right) }}{\left( \left| \lambda \right| -1\right)
!\;m\left( \lambda \right) !}J_{\lambda }\text{,}  \tag{3.16}
\end{equation}
with 
\begin{equation}
J_{\lambda }\left( q\right) =\left( \left| \lambda \right| -1\right)
!M_{\lambda }\left( q\right) q^{-n\left( \lambda \right) }\text{.} 
\tag{3.17}
\end{equation}
$J_{\lambda }$ is a polynomial with coefficients in $\mathbb{Z}$, of order
equal to zero and of degree equal to $\binom{\left| \lambda \right| -1}{2}%
+l\left( \lambda \right) -1-n\left( \lambda \right) $.
\end{corollary}

\begin{proof}
For $\lambda =\left( n\right) $, $\left( 3.16\right) $ gives $\left(
1.4\right) $ with $J_{\left( n\right) }=J_{n}$ which shows that $\left(
3.16\right) $ is a generalization of $\left( 1.4\right) $. Equations $\left(
3.16\right) $ and $\left( 3.17\right) $ follow easily from $\left(
3.1\right) $ and $\left( 3.2\right) $. The nullity of the order and the
value of the degree of $J_{\lambda }$ come from $\left( 3.5\right) $ and $%
\left( 3.3\right) $ respectively. Substituting the expression of $M_{\lambda
}$ coming from $\left( 3.17\right) $ into $\left( 3.10\right) $ gives for
all partitions $\lambda =\left( \lambda _{1},\lambda _{2},...,\lambda
_{r}\right) $ with $r\geq 2$%
\begin{equation}
J_{\lambda }=\left( \left| \lambda \right| -1\right) \binom{\left| \lambda
\right| -2}{\lambda _{r}-1}q^{-\left( r-1\right) \lambda _{r}}J_{\lambda
_{r}}J_{\lambda ^{\ast }}+\left( q-1\right)
\sum\limits_{i=1}^{r-1}q^{\lambda _{r}\left( i-r\right) }J_{\lambda ^{\left(
i\right) }}\text{,}  \tag{3.18}
\end{equation}
whith $\lambda ^{\ast }$ and $\lambda ^{\left( i\right) }$ defined above.
Equation $\left( 3.18\right) $ is a recurrence whose coefficients are
elements of $\mathbb{Z}\left[ q\right] $ possibly multipied by a negative
power of $q$. Note that the monomials with negative power of $q$ on the
right side of $\left( 3.18\right) $ necessarily cancel since the order of $%
J_{\lambda }$ is zero. Equation $\left( 3.18\right) $ makes it possible to
calculate all the polynomials $J_{\lambda }$ from the polynomials $J_{n}$
which belong to $\mathbb{N}\left[ q\right] $. So $J_{\lambda }$ are elements
of $\mathbb{Z}\left[ q\right] $.
\end{proof}

Equations $\left( 3.13\right) $ and $\left( 3.15\right) $ give directly with
Corolllary $\left( 3.5\right) $ :

\begin{corollary}
$J_{\lambda }$ polynomials are linked to $J_{n}^{\left( r\right) }$
polynomials by the following equations with strictly positive coefficients: 
\begin{equation}
\text{For \ }n\geq r\geq 1\text{,}\;\;\;\dfrac{\left( n-1\right) !}{r!\left(
n-r\right) !}q^{\binom{r}{2}}J_{n}^{\left( r\right) }=\sum\limits_{\left|
\lambda \right| =n,\;l\left( \lambda \right) =r}\dfrac{J_{\lambda }}{m\left(
\lambda \right) !}q^{n\left( \lambda \right) }\text{.}  \tag{3.19}
\end{equation}
\begin{equation}
\text{For }n\geq 2\text{,}\;\;\;\;J_{n}=\left( n-1\right)
\sum\limits_{\left| \lambda \right| =n-1}\dfrac{J_{\lambda }}{m\left(
\lambda \right) !}q^{n\left( \lambda \right) }\text{.}  \tag{3.20}
\end{equation}
\end{corollary}

With $\left( 3.19\right) $ it is possible to calculate the following
particular cases$\,$:$\ $%
\begin{equation}
\text{For all }n\geq 1\text{,}\;\;\;J_{1^{n}}=\left( n-1\right)
!J_{n}^{\left( n\right) }=\left( n-1\right) !\text{.}  \tag{3.21}
\end{equation}
\begin{equation}
\text{For all }n\geq 2\text{,}\;\;\;J_{1^{n-2}2}=\left( n-2\right)
!J_{n}^{\left( n-1\right) }=\left( n-2\right) !\left[ n-1\right] _{q}\text{.}
\tag{3.22}
\end{equation}

\begin{corollary}
\bigskip All the $J_{\lambda }$ can be recursively calculated from $J_{1}=1$
using Equations $\left( 3.18\right) $ and $\left( 3.20\right) .$
\end{corollary}

\begin{proof}
\bigskip By induction. Suppose we have calculated $J_{\lambda }$ for all $%
\lambda $ such that $\left| \lambda \right| \leq n-1$. Then, $J_{n}$ is
deduced from $\left( 3.20\right) $. For $\lambda $ such that $\left| \lambda
\right| =n$ with $l\left( \lambda \right) \geq 2$, we calculate them
successively in reverse lexicographic order using $\left( 3.18\right) $.
\end{proof}

We give below $J_{\lambda }$ for $\left| \lambda \right| \leq 4$ (Note that $%
J_{n}$ have already been given in $\left[ 8\right] $)

$\left| \lambda \right| =1\;\;J_{1}=1$.

$\left| \lambda \right| =2\;\;J_{2}=1$, $J_{1^{2}}=1$.

$\left| \lambda \right| =3\;\;J_{3}=2+q$, $J_{\left( 2,1\right) }=1+q$, \ $%
J_{1^{3}}=2$.

$\left| \lambda \right| =4$ \ $J_{4}=6+6q+3q^{2}+q^{3}$, $J_{\left(
3,1\right) }=3+3q+2q^{2}+q^{3}$, $J_{\left( 2,2\right) }=3+2q+q^{2}$, $%
J_{\left( 2,1,1\right) =1^{2}2^{1}}=2+2q+2q^{2}$, $J_{1^{4}}=6$.

\bigskip

\section{\protect\bigskip Conjectures about $J_{\protect\lambda }$}

Three heuristic arguments lead us to state the following conjectures.

\begin{conjecture}
\textit{For any partition }$\lambda $\textit{\ the coefficients of }$%
J_{\lambda }$\textit{\ are strictly positive.}
\end{conjecture}

To state the second and third conjectures, let us recall some definitions.
If $a=\left( a_{k}\right) _{0\leq k\leq n}$ is a finite sequence of real
numbers, then:\bigskip

* $a$ is unimodal if there is an index $j$, $0\leq j\leq n$ such that 
\begin{equation*}
a_{0}\leq a_{1}\leq ...\leq a_{j-1}\leq a_{j}\geq a_{j+1}\geq ...\geq a_{n}%
\text{.}
\end{equation*}

* $a$ is log-concave (resp. strictly log-concave) if 
\begin{equation*}
a_{j}^{2}\geq a_{j-1}a_{j+1}\ (\text{resp.}\;a_{j}^{2}>a_{j-1}a_{j+1})\ \ 
\text{for all}\ 1\leq j\leq n-1\text{.}
\end{equation*}
\qquad \qquad \qquad

* $a$ has no internal zeros if there do not exist $i<j<k$ satisfying $%
a_{i}\neq 0$, $a_{j}=0$ and $a_{k}\neq 0$.\bigskip

We will also say that the polynomial $P_{a}\left( q\right)
=a_{0}+a_{1}q+...+a_{n}q^{n}$ has one of these properties if the sequence of
its coefficients owns it. More details on theses properties and a rich
variety of unimodal and log-concave sequences can be found in $\left[ 5\text{%
, }7\text{, }17\right] $.

\begin{conjecture}
\textit{For any partition }$\lambda $\textit{, }$J_{\lambda }$\textit{\ is
log-concave.}
\end{conjecture}

If these two conjectures are true then $J_{\lambda }$ is also unimodal (see
Lemma 7.1.1 of $\left[ 5\right] $), thus we can also state:

\begin{conjecture}
For any partition $\lambda $, $J_{\lambda }$ is unimodal.
\end{conjecture}

Let us note that it is equivalent to formulate these conjectures for the
polynomials $J_{\lambda }$, $M_{\lambda }$ or $\widehat{J}_{\lambda }$=$%
J_{\lambda }/\Delta $, where $\Delta $ is the GCD of the coefficients of $%
J_{\lambda }$. The first argument supporting the statement of these
conjectures comes from the computations already performed. On the one hand,
the conjectures are obviously true for $J_{1^{n}}$ and $J_{1^{n-2}2}$ given
by Equations $\left( 3.21\right) $ and $\left( 3.22\right) $. On the other
hand, by developing a program based on Corollary $3.7$, Dr. Philippe Bodart
of the University of Lille, computed the polynomials $J_{\lambda }$ and
verified the conjectures for $\left| \lambda \right| =n\leq 43$. Given the
exponential growth in the number of integer partitions of $n$, this
represents 376325 polynomials. In other words, for each of these
polynomials, the strict positivity and log-concavity of their coefficients
were verified. Note that the highest degree among these polynomials is 861,
which correspond to the cases of $J_{43}$ and $J_{\left( 42,1\right) }$, as
one can calculate using Corollary 3.5.

The second argument in favor of these conjectures comes from the fact that
the coefficients of the polynomials $J_{n}^{\left( r\right) }$ are strictly
positive and log-concave, as we will demonstrate below. This follow from
Huh's work and is actually true for a broader class of polynomials denoted $%
I_{m}^{\left( a,b\right) }$ in $\left[ 20\right] $ ( Do not confuse the $a$
of $I_{m}^{\left( a,b\right) }$ with the sequence $a=\left( a_{k}\right)
_{0\leq k\leq n}$ above). We saw in ($\left[ 8\right] $, Section 6) that
these polynomials are related to $J_{n}^{\left( r\right) }$ through the
relation $J_{n}^{\left( r\right) }=I_{n-r}^{\left( r,1\right) }$. The
polynomials \ $I_{m}^{\left( a,b\right) }$ and their reciprocals, which are
the sum enumerators of \ the generalized classical parking functions, have
been the subject of extensive research, a summary of which can be found in
[20].

\begin{proposition}
\bigskip\ i) The sequence of the coefficients of each polynomial $%
I_{m}^{\left( a,b\right) }$ is strictly positive and log-concave, thus
unimodal.

\hspace{0.92in}ii) In the particular case $J_{n}^{\left( r\right)
}=I_{n-r}^{\left( r,1\right) }$, we have: 
\begin{equation}
J_{n}^{\left( r\right) }\left( q\right) =T_{n/r}\left( 1,q\right) \text{,} 
\tag{4.1}
\end{equation}
where $T_{n/r}$ is the Tutte polynomial of a graph which we denote $K_{n/r}$%
: $K_{n/r}$ is the complete graph $K_{n}$ in which all edges between $r$
given vertices are contracted.
\end{proposition}

\begin{proof}
\textit{i)} Let 
\begin{equation*}
I_{m}^{\left( a,b\right) }\left( q\right) =\sum\limits_{i=c}^{d}a_{i}q^{i}%
\text{.}
\end{equation*}
From the various properties of the polynomials $I_{m}^{\left( a,b\right) }$
and their reciprocals (see $\left[ 20\right] $ ), it is easy to see that $c=0
$, $d=ma+b\binom{m}{2}$ and 
\begin{equation}
a_{0}=m!\;\;\;a_{d}=1\text{.}  \tag{4.2}
\end{equation}
In [19] p. 662, it is shown \ that 
\begin{equation}
I_{m}^{\left( a,b\right) }\left( 1+t\right) =\sum\limits_{G^{\prime
}}t^{e\left( G^{\prime }\right) -m}\text{,}  \tag{4.3}
\end{equation}
where the sum is over the multicolor and connected graphs $G^{^{\prime }}$
without loop and whose vertex set is $V=\left\{ 0,1,2,...,m\right\} $. The
edges of $G^{\prime }$\ between two vertices $i,j\neq 0$ are to be taken
among the $b$ colored edges $\overline{0}$, $\overline{1}$, ..., $\overline{%
b-1}$ and those between $0$ and $i\neq 0$ are to be taken among the $a$
colored edges $\overline{0}$, $\overline{1}$,..., $\overline{a-1}$. $e\left(
G^{\prime }\right) $ is the number of edges of $G^{\prime }$. Let us now
consider the graph introduced in [16] p. 3115 with the notation $%
K_{m+1}^{b,a}$. This graph is defined as a complete graph on vertices $%
V=\left\{ 0,1,2,...,m\right\} $ with the edges $\left( i,j\right) $, $%
i,j\neq 0$ of multiplicity $b$ and the edges $\left( 0,i\right) $, $i\neq 0$
of multiplicity $a$. It is clear that there is a bijection between the
graphs $G^{\prime }$ and the connected spanning graphs of $K_{m+1}^{b,a}$.
The Tutte polynomial of $K_{m+1}^{b,a}$ is 
\begin{equation}
T\left( x,y\right) =\sum\limits_{A}\left( x-1\right) ^{c\left( A\right)
-1}\left( y-1\right) ^{c\left( A\right) +e\left( A\right) -(m+1)}\text{,} 
\tag{4.4}
\end{equation}
where the sum is over the spanning graphs $A$ of $K_{m+1}^{b,a}$, $c\left(
A\right) $ being the number of connected components of $A$, and $e\left(
A\right) $ its number of edges. Comparison of $\left( 4.3\right) $ and $%
\left( 4.4\right) $ shows that 
\begin{equation}
I_{m}^{\left( a,b\right) }\left( q\right) =T\left( 1,q\right) \text{.} 
\tag{4.5}
\end{equation}
Let us consider now the matroid $\mathcal{M}$ associated to $K_{m+1}^{b,a}$.
This matroid is graphical and therefore representable over any field (
Proposition 5.1.2 of $\left[ 15\right] $), so its dual $\mathcal{M}^{\ast }$%
also ( Corollary 2.2.9 of $\left[ 15\right] $). By noting $T^{\ast }$ the
Tutte polynomial of $\mathcal{M}^{\ast }$, we have $T^{\ast }\left(
y,x\right) =T(x,y)$, therefore $I_{m}^{\left( a,b\right) }\left( q\right)
=T^{\ast }\left( q,1\right) $. Let $\rho $ be the rank of $\ \mathcal{M}%
^{\ast }$ and $\left( h_{0},h_{1},...,h\rho \right) $ be the $h$-vector of
the matroid complex $IN\left( \mathcal{M}^{\ast }\right) $ (using Bj\"{o}%
rner's notation in $\left[ 3\right] $). It \ is known that ( see $\left[ 1%
\right] $ p.142): 
\begin{equation*}
\sum\limits_{k=0}^{\rho }h_{k}q^{\rho -k}=T^{\ast }\left( q,1\right) \text{.}
\end{equation*}
From $\left[ 11\text{, Theorem 3}\right] $, it follows that $I_{m}^{\left(
a,b\right) }$ is log-concave. Moreover, by the same theorem, we know that
its coefficients have non internal zeros. But these coefficients are
positive or zero (this comes for example, from their combinatorial
definition). The internal coefficients are therefore strictly positive, so
with $\left( 4.2\right) $ all the coefficients are strictly positive.
Unimodality follows from Lemma $7.1.1$ of [5].

\textit{ii) }First, note that the contracted graph $K_{n/r}$ clearly does
not depend on the choice of $r$ vertices (see $\left[ 9\right] $ for an
illustration). Let $K_{n}$ be the complete graph on $n$ vertices, and let $R$
be any set of $r$ vertices among these $n$ vertices. Contract all the edges
between the $r$ vertices of $R$, and number $0$ the vertex that replaces
these $r$ vertices. The set of other vertices $\mathbf{n}-R$, and the edges
between these vertices, remain unchanged in the contraction. In the
contracted graph, denoted $K_{n/r}$, there is thus one edge between two
vertices of $\mathbf{n}-R$, and $r$ edges between a vertex of $\mathbf{n}-R$
and vertex $0$. This shows that $K_{n/r}$ is isomorphic to $K_{n-r+1}^{1,r}$
and proves Equation $\left( 4.1\right) $ as a special case of $\left(
4.5\right) $.
\end{proof}

\begin{remark}
$K_{n}^{\left( r\right) }$ would be a more coherent notation to designate
the contracted graph of $ii)$ but this notation is already used for the
complete $r$-uniform hypergraph of order $n$, see [4, p. 310]. To our
knowledge, there is no simple standard notation assigned to the considered
family of graph. The proposed notation $K_{n/r}$ is quite suggestive, even
if it is not entirely satisfactory because a contraction is normally done on
edges and not on vertices.
\end{remark}

We have just proven that the polynomials $J_{n}^{\left( r\right) }$, and in
particular $J_{n}=J_{n}^{\left( 1\right) }$, verify the conjectures.
Moreover, there are close links between the polynomials $J_{n}^{\left(
r\right) }$ and $J_{\lambda }$. Firstly, they have a non-empy intersection
corresponding to the three sequences of polynomials extracted from $%
J_{\lambda }$, i.e., for $\lambda =\left( n\right) $ as seen before, and for 
$\lambda =$ $1^{n}$ and $1^{n-2}2^{1}$ as shown by Equations $\left(
3.21\right) $ and $\left( 3.22\right) $, respectively. Furthermore,
Equations $\left( 3.19\right) $ and $\left( 3.20\right) $ are equations
relating the $J_{\lambda }$ to the $J_{n}^{\left( r\right) }$ with
coefficients that are polynomials in $q$ with strictly positive integer
coefficients. All these facts constitute the second argument in favor of the
conjectures. We now present the third argument.

\section{\protect\bigskip Pascalian part of $M_{\protect\lambda }$}

\begin{theorem}
Let $\lambda $ be a partition of the integer $n$ with $l\left( \lambda
\right) =r$ \ and let 
\begin{equation}
v\left( \lambda \right) =d\left( \lambda \right) -n+2=\binom{n-2}{2}+r-1%
\text{.}  \tag{5.1}
\end{equation}
For $\lambda \neq 1^{n}\Leftrightarrow 1\leq r\leq n-1$, the part of $%
M_{\lambda }$ defined in Theorem 3.1, of degree $\geq $ $v\left( \lambda
\right) $, is given by the following polynomial which only depends on $n$
and $r$: 
\begin{equation}
P_{n}^{\left( r\right) }\left( q\right) =\dfrac{\left( r-1\right) !}{\left(
n-1\right) !}\sum\limits_{i=0}^{n-2}\binom{n-r-1+i}{n-r-1}q^{d\left( \lambda
\right) -i}=\dfrac{\left( r-1\right) !}{\left( n-1\right) !}%
\sum\limits_{j=0}^{n-2}\binom{2n-r-3-j}{n-2-j}q^{v\left( \lambda \right) +j}.
\tag{5.2}
\end{equation}
\end{theorem}

$v\left( \lambda \right) $ is the order of $P_{n}^{\left( r\right) }$.
Theorem 5.1 shows that the $n-1$ coefficients of $\left( \left( n-1\right)
!/(r-1)!\right) M_{\lambda }$ with highest degree are given by the first $%
n-1 $ coefficients of column $\left( n-r-1\right) $\ in Pascal's triangle. $%
P_{n}^{\left( r\right) }$ and it coefficients will be called Pascalian part
and Pascalian coefficients of $M_{\lambda }$ (likewise for the parts and
coefficients corresponding to $J_{\lambda }$ and $\widehat{J}_{\lambda }$).

\begin{proof}
We prove it by a double induction on $n$ and $r$.

a) If $n=2$ and $r=1$, we clearly have $M_{\left( 2\right) }=J_{2}=1$ and $%
n\left( \left( 2\right) \right) =d\left( \left( 2\right) \right) =0$. Thus
the assertion of Theorem 5.1 is true for $n=$ $\left| \lambda \right| =2$
and $\lambda \neq 1^{2}$.

b) Suppose that for $k\leq n$, \ $\left| \lambda \right| =k$ and $\lambda
\neq 1^{k}$ the assertion of Theorem 5.1 is true. For $\lambda =\left(
n+1\right) $, $\left( 3.15\right) $ gives 
\begin{equation}
M_{\left( n+1\right) }=\sum\limits_{r=1}^{n}\sum\limits_{\left| \lambda
\right| =n\;l\left( \lambda \right) =r}\dfrac{M_{\lambda }}{m\left( \lambda
\right) !}\text{.}  \tag{5.3}
\end{equation}
For $\lambda $ such that $\left| \lambda \right| =n$ and different from $%
1^{n}$, we have $r\leq n-1$, so $\left( 5.1\right) $ gives 
\begin{equation*}
v\left( \lambda \right) =d\left( \lambda \right) -n+2=\binom{n-1}{2}%
+r-\left( n-1\right) \leq \binom{n-1}{2}=v\left( \left( n+1\right) \right) 
\text{.}
\end{equation*}
Therefore, in the sum deduced from $\left( 5.3\right) $\ which gives the
coefficient of $q^{m}$ in $M_{\left( n+1\right) }$,\ all coefficients $%
\left\langle q^{m}\right\rangle M_{\lambda }$ for $\lambda \neq 1^{n}$ are
Pascalian provided that $m\geq v\left( \left( n+1\right) \right) $. We saw
that 
\begin{equation*}
\dfrac{M_{1^{n}}}{n!}=\dfrac{\widetilde{m}_{1^{n}}}{n!}=\dfrac{q^{\binom{n}{2%
}}}{n!}\text{,}
\end{equation*}
and according to $\left( 3.3\right) $, $d\left( \left( n+1\right) \right) =%
\binom{n}{2}$. Thus for $v\left( \left( n+1\right) \right) \leq m\leq
d\left( \left( n+1\right) \right) -1$, the terms of the sum in $\left(
5.3\right) $ which contribute to $\left\langle q^{m}\right\rangle M_{\left(
n+1\right) }$ are the terms associated with the partitions verifying $%
d\left( \lambda \right) \geq m$, i.e. with $\left( 3.3\right) $, $r=l\left(
\lambda \right) \geq m+1-\binom{n-1}{2}$. By $\left( 5.2\right) $ and $%
d\left( \lambda \right) -i=m\Leftrightarrow i=d\left( \lambda \right) -m$,
it follows that 
\begin{equation*}
\left\langle q^{m}\right\rangle M_{\lambda }=\dfrac{\left( r-1\right) !}{%
\left( n-1\right) !}\binom{n-r-1+d\left( \lambda \right) -m}{d(\lambda )-m}%
\text{,}
\end{equation*}
by adding these expressions, this gives 
\begin{equation*}
\left\langle q^{m}\right\rangle M_{\left( n+1\right) }=\sum\limits_{r=m+1-%
\binom{n-1}{2}}^{n}\dfrac{\left( r-1\right) !}{\left( n-1\right) !}%
\sum\limits_{\left| \lambda \right| =n,l\left( \lambda \right) =r}\dfrac{%
\binom{n-r-1+d\left( \lambda \right) -m}{d(\lambda )-m}}{m\left( \lambda
\right) !}\text{.}
\end{equation*}
Let $m=v\left( \left( n+1\right) \right) +j$, substituting $m$ in the above
equation, we obtain for $0\leq j\leq n-2$: 
\begin{equation}
\left\langle q^{v\left( \left( n+1\right) \right) +j}\right\rangle M_{\left(
n+1\right) }=\sum\limits_{r=j+1}^{n}\dfrac{\left( r-1\right) !}{\left(
n-1\right) !}\binom{n-2+j}{r-1-j}\sum\limits_{\left| \lambda \right|
=n,l\left( \lambda \right) =r}\dfrac{1}{m\left( \lambda \right) !}\text{.} 
\tag{5.4}
\end{equation}
Furthermore, it is known that 
\begin{equation}
\sum\limits_{\left| \lambda \right| =n,l\left( \lambda \right) =r}\dfrac{1}{%
m\left( \lambda \right) !}=\left( n!\right) ^{-1}B_{n,r}\left(
1!,2!,...\right) =\binom{n-1}{r-1}\dfrac{1}{r!}\text{,}  \tag{5.5}
\end{equation}
where $B_{n,r}\left( x_{1},x_{2},...\right) $ is the Bell's partial
exponential polynomial, the last equality coming from equation $\left[ 3h%
\right] $ p.135 of $\left[ 10\right] $. With $\left( 5.5\right) $ and the
change of index $l=r-1-j$, $\left( 5.4\right) $ successively becomes 
\begin{equation*}
\left\langle q^{v\left( \left( n+1\right) \right) +j}\right\rangle M_{\left(
n+1\right) }=\dfrac{1}{n!}\sum\limits_{r=j+1}^{n}\binom{n-2-j}{r-1-j}\binom{n%
}{r}=\dfrac{1}{n!}\sum\limits_{l=0}^{n-1-j}\binom{n-2-j}{l}\binom{n}{n-1-j-l}%
\text{.}
\end{equation*}
Using\ Chu-Vandermonde's identity given by (see Formula $\left[ 13c^{\prime }%
\right] $ p. 44 of $\left[ 10\right] $): 
\begin{equation*}
\sum\limits_{l=0}^{L}\binom{A}{l}\binom{B}{L-l}=\binom{A+B}{L}\text{,}
\end{equation*}
we finally obtain for $0\leq j\leq n-2$ 
\begin{equation*}
\left\langle q^{v\left( \left( n+1\right) \right) +j}\right\rangle M_{\left(
n+1\right) }=\dfrac{1}{n!}\binom{2n-2-j}{n-1-j}
\end{equation*}
which is indeed the coefficient of the last member of $\left( 5.2\right) $
for $\lambda =\left( n+1\right) $. In the case $m=d\left( \left( n+1\right)
\right) =\binom{n}{2}\Leftrightarrow j=n-1$ we already know by $\left(
3.4\right) $ that $\left\langle q^{d\left( \left( n+1\right) \right)
}\right\rangle M_{\left( n+1\right) }=1/n!$ which is the value given by $%
\left( 5.2\right) $. We have therefore proved Theorem 5.1 for $\lambda
=\left( n+1\right) $.

c) Let us now make an induction on $r=l\left( \lambda \right) $ with $\left|
\lambda \right| =n+1$. The case $r=1$ is proved by b). Assume that the
assertion of Theorem 5.1 is true for all $l$ such that $1\leq l\leq r\leq n-1
$ and let show that it is true for $r+1.$ If $\lambda =\left( \lambda
_{1},\lambda _{2},...,\lambda _{r},\lambda _{r+1}\right) $ with $\left|
\lambda \right| =n+1$ and $r+1\leq n$, $\left( 3.10\right) $ gives 
\begin{equation}
M_{\lambda }=M_{\left( \lambda _{r+1}\right) }M_{\lambda ^{\ast }}+\left(
q-1\right) \sum\limits_{i=1}^{r}M_{\lambda ^{\left( i\right) }}\text{.} 
\tag{5.6}
\end{equation}
According to Theorem 3.1, the degree of $M_{\left( \lambda _{r+1}\right)
}M_{\lambda ^{\ast }}$ is 
\begin{equation*}
\delta =d\left( \left( \lambda _{r+1}\right) \right) +d\left( \lambda
_{1},\lambda _{2},...,\lambda _{r}\right) =\binom{\lambda _{r+1}-1}{2}+%
\binom{n+1-\lambda _{r+1}-1}{2}+r-1\text{,}
\end{equation*}
and $\left( 5.1\right) $ gives here 
\begin{equation}
v\left( \lambda \right) =d\left( \lambda \right) -\left| \lambda \right| +2=%
\binom{n-1}{2}+r\text{.}  \tag{5.7}
\end{equation}
An easy calculation shows that $v\left( \lambda \right) -\delta =\left|
\lambda \right| \left( \lambda _{r+1}-1\right) -\lambda _{r+1}^{2}+2$. As $%
\left| \lambda \right| \geq \lambda _{r+1}+1$, it follows $v\left( \lambda
\right) -\delta \geq 1$, so $v\left( \lambda \right) >\delta $. Therefore,
in the Pascalian part of $M_{\lambda }$ only the term $\left( q-1\right)
\sum\limits_{i=1}^{r}M_{\lambda ^{\left( i\right) }}$ of $\left( 5.6\right) $
intervenes. For $1\leq i\leq r$ we have\ $l\left( \lambda ^{\left( i\right)
}\right) =r$ and $\left| \lambda ^{\left( i\right) }\right| =n+1$, so the $%
\lambda ^{\left( i\right) }$ verify the assertion of Theorem 5.1 and their
Pascalian parts are equal to $P_{n+1}^{\left( r\right) }$ with 
\begin{equation*}
v\left( \lambda ^{\left( i\right) }\right) =\binom{n-1}{2}+r-1\text{,}
\end{equation*}
which implies with $\left( 5.7\right) $ 
\begin{equation*}
v\left( \lambda \right) =v\left( \lambda ^{\left( i\right) }\right) +1\text{.%
}
\end{equation*}
Hence, for $m$ such that $v\left( \lambda \right) \leq m\leq d\left( \lambda
\right) -1$, we have 
\begin{equation*}
\left\langle q^{m}\right\rangle M_{\lambda }\left( q\right) =\left\langle
q^{m}\right\rangle \left( q-1\right) rP_{n+1}^{\left( r\right) }\left(
q\right) \text{,}
\end{equation*}
i.e. with $\left( 5.2\right) $ 
\begin{equation}
\left\langle q^{m}\right\rangle M_{\lambda }\left( q\right) =\dfrac{r!}{n!}%
\left\langle q^{m}\right\rangle \left( q-1\right) \sum\limits_{j=0}^{n-1}%
\binom{2n-r-1-j}{n-1-j}q^{v\left( \lambda ^{\left( i\right) }\right) +j}%
\text{.}  \tag{5.8}
\end{equation}
In the right hand side of $\left( 5.8\right) $, the polynomial after $%
\left\langle q^{m}\right\rangle $ is written by expanding: 
\begin{eqnarray*}
&&-\binom{2n-r-1}{n-1}q^{v\left( \lambda ^{\left( i\right) }\right)
}+\sum_{j=1}^{n-1}\left[ \binom{2n-r-1-\left( j-1\right) }{n-1-\left(
j-1\right) }-\binom{2n-r-1-j}{n-1-j}\right] q^{v\left( \lambda ^{\left(
i\right) }\right) +j} \\
&&+\binom{2n-r-1-\left( n-1\right) }{n-1-\left( n-1\right) }q^{v\left(
\lambda ^{\left( i\right) }\right) +n}\text{.}
\end{eqnarray*}
With the classical relation $\binom{N}{p}=\binom{N-1}{p-1}+\binom{N-1}{p}$,
the square bracket in the above equation becomes: 
\begin{equation*}
\binom{2\left( n+1\right) -\left( r+1\right) -3-\left( j-1\right) }{%
n+1-2-\left( j-1\right) }\text{.}
\end{equation*}
By making the change of index $j^{\prime }=j-1$ and keeping only the
monomials of $M_{\lambda }$\ with degreee $\geq v\left( \lambda \right) $,
we obtain for the sum of these monomials 
\begin{equation*}
\dfrac{r!}{n!}\sum\limits_{j^{\prime }=0}^{n-1}\binom{2\left( n+1\right)
-\left( r+1\right) -3-j^{\prime }}{\left( n+1\right) -2-j^{\prime }}%
q^{v\left( \lambda \right) +j^{\prime }}\text{,}
\end{equation*}
which is indeed the Pascalian part $P_{n+1}^{\left( r+1\right) }$ given by $%
\left( 5.2\right) $.
\end{proof}

Here is now the proof of the conjectures of Section 4 for part of the
coefficients of $J_{\lambda }$, which constitutes the third announced
argument.

\begin{corollary}
For any partition $\lambda $ such that $\lambda \neq 1^{\left| \lambda
\right| }$, the sequence of the last $\left| \lambda \right| -1$
coefficients of $J_{\lambda }$ is strictly positive and log-concave, then
unimodal. If $l\left( \lambda \right) <\left| \lambda \right| -1$ this
sequence is in fact strictly log-concave.
\end{corollary}

\begin{proof}
It is known that the coefficients of the columns\ of Pascal's triangle are
strictly positif and log-concave. More precisely for any pair of integers $%
\left( N,p\right) $ such that $p\geq 0$ et $N\geq p+1$%
\begin{equation*}
\binom{N}{p}^{2}-\binom{N+1}{p}\binom{N-1}{p}=\left( \dfrac{N!}{p!\left(
N-p\right) !}\right) ^{2}p\text{.}
\end{equation*}

Thus, the columns of Pascal's triangle are strictly log-concaves for $p\geq 1
$, the column $p=0$, consisting of 1, is only log-concave. The corollary
follows for $J_{\lambda }$ and $\lambda \neq 1^{\left| \lambda \right| }$
from the proportionality, given by Theorem 5.1 and $\left( 3.17\right) $,
between Pascalian coefficients of $J_{\lambda }$ and those of column $\left(
n-r-1\right) $\ of Pascal's triangle.
\end{proof}

Incidentally it is easy to obtain a result similar to Theorem 5.1 for the
polynomials $J_{n}^{\left( r\right) }$:

\begin{corollary}
For $1\leq r\leq n-1$, and $0\leq i\leq n-2$, we have: 
\begin{equation}
\left\langle q^{\binom{n-1}{2}-\binom{r-1}{2}-i}\right\rangle J_{n}^{\left(
r\right) }=\binom{n-r-1+i}{n-r-1}\text{,}  \tag{5.9}
\end{equation}
which means that the $\left( n-1\right) $ coefficients of $J_{n}^{\left(
r\right) }$ of the highest degree are given by the first coefficients of
column $\left( n-r-1\right) $\ of Pascal's triangle.
\end{corollary}

\begin{proof}
\bigskip Equations $\left( 3.13\right) $ and $\left( 5.2\right) $
successively give for $n-1\geq r\geq 1$, $m=d\left( \lambda \right) -i$ and $%
0\leq i\leq n-2$: 
\begin{equation*}
\left\langle q^{m}\right\rangle q^{\binom{r}{2}}\dfrac{J_{n}^{\left(
r\right) }}{r!\left( n-r\right) !}=\sum\limits_{\left| \lambda \right|
=n,\;l\left( \lambda \right) =r}\dfrac{\left\langle q^{m}\right\rangle
M_{\lambda }}{m\left( \lambda \right) !}=\dfrac{\left( r-1\right) !}{\left(
n-1\right) !}\binom{n-r-1+i}{n-r-1}\sum\limits_{\left| \lambda \right|
=n,\;l\left( \lambda \right) =r}\dfrac{1}{m\left( \lambda \right) !}\text{.}
\end{equation*}

With $\left( 5.5\right) $ and after simplification we obtain: 
\begin{equation*}
\left\langle q^{m}\right\rangle q^{\binom{r}{2}}J_{n}^{\left( r\right)
}=\left\langle q^{m-\binom{r}{2}}\right\rangle J_{n}^{\left( r\right) }=%
\binom{n-r-1+i}{n-r-1}\text{,}
\end{equation*}

and $\left( 3.3\right) $ implies 
\begin{equation*}
m-\binom{r}{2}=\binom{n-1}{2}+r-1-i-\binom{r}{2}=\binom{n-1}{2}-\binom{r-1}{2%
}-i\text{,}
\end{equation*}
which gives $\left( 5.9\right) $. Moreover, Theorem 6.3 of $\left[ 8\right] $
shows that $\binom{n-1}{2}-\binom{r-1}{2}$ is the degree of $J_{n}^{\left(
r\right) }$, which ends the proof.
\end{proof}

\section{\protect\bigskip To go further}

The numerical calculations we performed were implemented on a personal
computer. The value $n=\left| \lambda \right| =43$ represents the limit that
can be reached within a reasonable computation time. With the three
arguments presented, the author has little doubt about the validity of the
conjectures (even though, of course, there is no certainty until they have
been proven). In any case, the programming codes are available upon request
for anyone wishing to use more powerful means to verify these conjectures
for a larger integer $n$.

The main task now is to prove Conjecture 4.1. Two approaches can be
considered for this purpose. The first consists in generalizing the
inductive proof from Section 5. Ideally, the goal would be to find a formula
that extends the one obtained for Pascalian coefficients so as to include
all coefficients. It would already be beneficial to discover a recurrence
relation between the $J_{\lambda }$ with coefficients in $\mathbb{N}\left[ q%
\right] $ - unlike Recurrence $\left( 3.18\right) $ which has coefficients
in $\mathbb{Z}\left[ q\right] $- in order to prove Conjecture 4.1. The
second approach, of a combinatorial nature, would consists in associating to
each integer partition $\lambda $, a mathematical object, let's call it $%
K_{\lambda }$, whose polynomial $J_{\lambda }$ (or $\widehat{J}_{\lambda }$%
)\ would be the enumerator of some quantity related to that object. Thus, $%
K_{\lambda }$ would be a generalization of the complete graph $K_{n}$
associated with $J_{n}$, but in a different direction from the contracted
graph $K_{n/r}$ associated with $J_{n}^{\left( r\right) }$. It is noteworthy
that it seems to emerge from what follows that $K_{\lambda }$ cannot be a
graph.

Suppose that we have discovered this object $K_{\lambda }$. A fairly natural
approach to prove Conjecture 4.2 would then be to reproduce the method used
in the proof of Proposition 4.4. This would involve associating to each $%
K_{\lambda }$ a matroid $M_{\lambda }$ whose Tutte polynomial would satisfy
the following analogue of $\left( 4.5\right) $: 
\begin{equation*}
\widehat{J}_{\lambda }\left( q\right) =T_{\lambda }\left( 1,q\right) \text{.}
\end{equation*}

Note, in this regard, that the log-concavity of the $h$-vector of the
independence complex of a matroid has recently been demonstrated for all
matroid [2]. This avoids being limited to graphic matroids for $\mathcal{M}%
_{\lambda }$ (and therefore to graphs for $K_{\lambda }$). However, this
approach presents a difficulty. Indeed, we have calculated that: 
\begin{equation*}
J_{\left( 3,2,1\right) }\left( q\right)
=10+30q+35q^{2}+35q^{3}+30q^{4}+20q^{5}+12q^{6}+6q^{7}+2q^{8}\text{.}
\end{equation*}

Since the coefficients $2$ and $35$ are coprime, $J_{\left( 3,2,1\right) }$
cannot be reduced to a monic polynomial with integer coefficients. But if we
were to find a matroid $\mathcal{M}$ whose Tutte polynomial $T$ satisfies $%
J_{\left( 3,2,1\right) }\left( q\right) =T\left( 1,q\right) $, then
according to [1, Eq. (1.23) p.138], we would have:

\begin{equation*}
T\left( 1,q\right) =\sum \left( q-1\right) ^{\left| A\right| -\rho }\text{,}
\end{equation*}
where the sum is taken over all the subsets $A$ of the ground set $E$\ of $%
\mathcal{M}$, whose ranks are equal to $\rho $. The highest-degree monomial
of $T\left( 1,q\right) $ is therefore clearly $q^{\left| E\right| -\rho }$,
so $T\left( 1,q\right) $ would be monic leading to a contradiction.

To put this difficulty into perspective, it is worth noting that the
techniques used to prove Huh's theorem, mentionned in Section 4, have been
considerably generalized (see, for example, [6]). It may therefore be
conceivable to overcome the difficulty encountered by employing these
techniques.

It is worth emphasizing, however, that before embarking on the proof of
Conjecture 4.2 on log-concavity, a preliminary step is, as mentionned
earlier, to prove Conjecture 4.1 on the strict positivity of the
coefficients of $J_{\lambda }$. The latter poses an interesting challenge
from a combinatorial point of view.

\textit{I am grateful to Philippe Bodart for the computation of the }$%
J_{\lambda }$\textit{\ polynomials and the verification of the conjectures
up to }$n=43$\textit{, using the programs he developped.}

\textit{\bigskip }

\textbf{References}

$\left[ 1\right] $ F. Ardila, Algebraic and Geometric Methods in Enumerative
Combinatorics. In \textit{Handbook of Enumerative Combinatorics}, CRC Press,
Boca Raton, FL (2015) pp. 3-158.

$\left[ 2\right] $ F. Ardila, G. Denham, and J. Huh., Lagrangien geometry of
matroids. \textit{J. Amer. Math. Soc.,} \textbf{36} (2023) 727-794.

$\left[ 3\right] $ A. Bj\"{o}rner, The Homology and Shellability of Matroids
and Geometric Lattices. In \textit{Matroid Applications}, N. White (ed),
Cambridge University Press. pp. 226-283 (1992).

$\left[ 4\right] $ J.A. Bondy, U.S.R Murty, \textit{Graph Theory}. Springer,
London (2008).

$\left[ 5\right] $ P. Br\"{a}nden, Unimodality, Log-concavity,
Real-rootedness and Beyond. In \textit{Handbook of Enumerative Combinatorics}%
, CRC Press, Boca Raton, FL (2015) pp. 437-483.

$\left[ 6\right] $ P. Br\"{a}nden, J. Huh, Lorentzian Polynomials. \textit{%
Ann. of Math}. (2), \textbf{192} (3), 821-891, (2020)

$\left[ 7\right] $ F. Brenti, Log-concave and unimodal sequences in algebra,
combinatorics, and geometry: an update. \textit{Comtemp. Math.}, vol. 
\textbf{178} (1994), pp. 71-89.

$\left[ 8\right] $ V. Brugidou, A q-analog of certain symmetric functions
and one of its specializations, (2023); \textbf{arXiv}:2302.11221.

$\left[ 9\right] $ V. Brugidou, On a particular specialization of monomial
symmetric functions, International Conference on Enumerative Combinatorics
and Applications (ICECA), (September 4-6-, 2023).

$\left[ 10\right] $ L. Comtet, \textit{Advanced Combinatorics}. Springer
Netherlands, Dordrecht,1974. $\ $

$\left[ 11\right] $ J. Huh, \textit{h-}vectors of matroids and logarithmic
concavity. \textit{Adv. Math.}, \textbf{270 }(2015) 49-59.

$\left[ 12\right] $ I. G. Macdonald, \textit{Symmetric functions and Hall
polynomials}, second ed., Oxford University Press, New York, 1995.

$\left[ 13\right] $ C. L. Mallows and J. Riordan, The inversion enumerator
for labeled trees. \textit{Bull. Amer. Soc. }\textbf{74} (1968) 92-94.

$\left[ 14\right] $ M. Merca, Augmented monomials in terms of power sums, 
\textit{SpringerPlus }(2015) 4:724 DOI 10.11686/s40064-015-1506-5.

$\left[ 15\right] $ J. Oxley, \textit{Matroid Theory}, second ed., Oxford
Graduate Texts in Mathematic, Oxford university Press, 2011.

$\left[ 16\right] $ A. Postnikov and B. Shapiro, Trees, parking functions,
syzygies, and deformations of monomial ideals. \textit{Trans. of the Amer.
Math. Soc.}, \textbf{356} Number 8 (2004) pp.3109-3142.

$\left[ 17\right] $ R.P. Stanley, Log-concave and unimodal sequences in
algebra, combinatorics, and geometry. \textit{Ann. New York Acad. Sci.} 
\textbf{576} (1989) 500-535.

$\left[ 18\right] $ R. P. Stanley, Hyperplane arrangements, parking
functions, and tree inversions. In \textit{Mathematical Essays in Honor of
Gian-carlo Rota (Cambridge, MA, 1996)}, vol. 161 of Progr. Math.,
Birh\"{a}user, Boston, MA (1998) 359-375.

$\left[ 19\right] $ C. H.Yan, Generalized parking functions, tree
inversions, and multicolored graphs. \textit{Adv. in Appl. Math}. \textbf{27}%
(2-3) (2001) 641-670

$\left[ 20\right] $ C. H. Yan, Parking Functions. In \textit{Handbook of
Enumerative Combinatorics}, CRC Press, Boca Raton, FL (2015) pp. 835-893.

\end{document}